\documentclass[twocolum,fleqn]{autart}    

\usepackage{graphicx}          

\usepackage{amssymb}
\usepackage[cmex10]{amsmath}
\usepackage{accents}
\usepackage{epsfig}

\usepackage{bm}
\usepackage{amsmath}
\usepackage{color}
\usepackage{dsfont}
\usepackage{empheq}
\usepackage{hyperref}
\usepackage{array}
\usepackage{cite}
\usepackage{subfig}
\usepackage{pgfplots}
\usepackage{tikz}
\usetikzlibrary{arrows,shapes,patterns,positioning,calc}
\tikzstyle{boxstyle}=[draw=black,inner sep=7pt]     
\tikzstyle{arrowstyle}=[]    
\usetikzlibrary{arrows, decorations.markings}

\newlength\figureheight
\newlength\figurewidth

\newtheorem{remark}{Remark}

\newcommand{\ba}{\begin{align}}
\newcommand{\ea}{\end{align}}
\newcommand{\fr}{\frac}

\newcommand{\ep}{\varepsilon}

\begin{document}

\begin{frontmatter}

\title{Arctic Sea Ice State Estimation \\From Thermodynamic PDE Model} 
 

\author[UCSD]{Shumon Koga}\ead{skoga@eng.ucsd.edu},               
\author[UCSD]{Miroslav Krstic}\ead{krstic@ucsd.edu}  

\address[UCSD]{Department of Mechanical and Aerospace Engineering, University of California, San Diego, La Jolla, CA 92093-0411 USA}

\begin{keyword}                           
Arctic sea ice, State estimation, Moving boundaries, Backstepping
\end{keyword}

\begin{abstract}                          
Recent rapid loss of the Arctic sea ice motivates the study of the Arctic sea ice thickness. Global climate model that describes the ice's thickness evolution requires an accurate spatial temperature profile of the Arctic sea ice. However, measuring the complete temperature profile is not feasible within and throughout the Arctic icecap. Instead, measuring the ice's thickness is doable with the acquisition of data from submarine and satellite devices. In this paper, we develop a backstepping observer algorithm to estimate the temperature profile for the Arctic sea ice model via available measurements of sea ice thickness and sea ice surface temperature. The observer is designed in a rigorous manner to drive the temperature profile estimation error to zero, for a salinity-free sea ice model. Moreover, the proposed  observer is  used to estimate the temperature profile of the original sea ice model with salinity via numerical simulation. In comparison with the straightforward open-loop algorithm, the simulation results illustrate that our observer design achieves ten times faster convergence of the estimated temperature. 
\end{abstract}

\end{frontmatter}

\section{Introduction}
\subsection{Motivation}
The Arctic sea ice has been studied intensively  in the field of climate and geoscience. One of the main reasons is due to ice-albedo feedback which influences climate dynamics through the high reflectivity of sea ice. The other reason is the rapid decline of the Arctic sea ice extent in the recent decade shown in several observations. These observations motivate the investigation of future sea ice amount. Several studies have developed a computational model of the Arctic sea ice and performed numerical simulations of the model with initial sea ice temperature profile. However, the spatially distributed temperature in sea ice is difficult to recover in realtime using a limited number of thermal sensors. Hence, the online estimation of the sea ice temperature profile based on some available measurements is crucial for the prediction of the sea ice thickness.

\subsection{Literature}
A thermodynamic model for the Arctic sea ice was firstly developed in \cite{maykut1971} (hereafter MU71), in which the authors investigated the correspondence between the annual cycle pattern acquired from the simulation and the empirical data of \cite{unter1961}. The model involves a temperature diffusion equation evolving on a spatial domain defined as the sea ice thickness. Due to melting or freezing phenomena, the aforementioned spatial domain is time-varying.  Such a model is called ``Stefan problem" \cite{Gupta03} which is described by a parabolic partial differential equation (PDE) with a state-dependent moving boundary driven by a Neumann boundary value.

Refined models of MU71 have been suggested in literature. For instance, \cite{Semtner1976} proposed a numerical model to achieve faster and accurate computation of MU71 by discretizing the temperature profile into some layers and neglecting salinity effect. The energy conserving model of MU71 was introduced in \cite{Bitz1999} by taking into account the internal brine pocket melting on surface ablation and the vertically varying salinity profile. Their thermodynamic model was demonstrated by \cite{Bitz2001} using  a global climate model with a  Lagrangian ice thickness distribution. Combining these two models, \cite{Winton2000} developed an energy conserving three-layer model of sea ice by treating the upper half of the ice as a variable heat capacity layer. 

Remote sensing techniques have been employed to obtain the Arctic sea ice data in several studies. In \cite{Hall2004}, authors suggested an algorithm to calculate sea ice surface temperature using the satellite measured brightness temperatures, which provided an excellent measurement of the actual surface temperature of the sea ice during the Arctic cold period. The Arctic sea ice thickness data were acquired in \cite{Kwok2009} through a satellite called "ICESat" during 2003-2008, and compared with the data in \cite{Rothrock2008} observed by submarine during 1958-2000. More recent data describing the evolution of the sea ice thickness have been collected between 2010 and 2014 from the satellite called "CryoSat-2" \cite{Laxon2013, Kwok2015}. 

On the other hand, state estimation has been studied as a specific type of data assimilation which utilizes the numerical model along with measured value. For finite dimensional systems associated with noisy uncertain measurement, a systematic well known approach is Kalman Filter. Another well known method is Luenberger type state observer which reconstructs the state variable from partially measured variables. For the application to sea ice, \cite{Fenty2013} developed an adjoint based method as an iterative state and parameter estimation for the coupled sea ice-ocean in the Labrador Sea and Baffin Bay to minimize an uncertainty-weighted model-data misfit in a least square sense as suggested in \cite{Wunsch1996, Wunsch2007} using Massachusetts Institute of Technology general circulation model (MITgcm) developed in \cite{Marshall1997}. In \cite{Fenty15}, the same methodology was applied to reconstruct the global ocean and ice concentration. Their sea ice model is based on the zero-layer approximation of numerical model in \cite{Semtner1976}, which is a crude model lacking an internal heat storage and promoting fast melting. 

State estimation for infinite dimensional systems is much more challenging than the estimation of finite dimensional systems. The infinite dimensional systems are often modeled by means of PDEs which describes a time evolution of spatially distributed state variables. A systematic estimation method for PDEs was developed via "backsteppping design" which employs the state observer structure for PDEs, see \cite{miroslav08,miroslav09,susto10,tang11,andrey05} for instance. Recently, the backstepping method for PDEs has been employed to the Stefan problem for the design of boundary controllers in \cite{Shumon16state} and estimators with boundary measurements in \cite{Shumon16output,Shumon17journal}. As further extensions, \cite{Shumon17ACC} designed a state feedback control for the Stefan problem under the material's convection, \cite{koga2017CDC,koga_2019delay} designed a delay compensated control for the Stefan problem, and \cite{koga2018ISS} investigated an input-to-state stability of the control of Stefan problem in \cite{Shumon16state} with respect to an unknown heat loss at the interface. 

\subsection{Contributions}
Our conference paper \cite{Shumon17seaice} applied the state estimation for the Stefan problem developed in \cite{Shumon17journal} to the thermodynamic model of Arctic sea ice in \cite{maykut1971}. The present paper improves the design of \cite{Shumon17seaice} by removing the availability of the temperature gradient at the ice-ocean interface which is nearly impossible to measure. The designed algorithm provides the temperature profile estimation in Arctic sea ice via backstepping method. The observer utilizes the available measurements of the sea ice thickness and the ice surface temperature. A simplified model for the Arctic sea ice is formulated by neglecting the salinity effect and  the exponential convergence of its  temperature estimation is ensured by using the designed observer. The simulation study for the proposed observer along with the original model is performed  in order to investigate the convergence performance numerically. The simulation results illustrate that the estimated temperature profile converges uniformly to the actual sea ice temperature ten times faster than the straightforward open-loop estimator. 
\subsection{Organization} 
This paper is organized as follows: Section \ref{sec:2} describes the thermodynamic model of the Arctic sea ice in MU71 and introduces its simplification. Section \ref{sec:3} develops the state estimation design for the simplified model via backstepping PDE observer and shows the exponential stability of the estimation error system. Section \ref{sec:4} illustrates the simulation of the designed observer applied to the original model of MU71. The paper ends with the concluding remarks and future direction in Section \ref{sec:5}.

\section{Thermodynamic Model of Sea Ice}\label{sec:2}
The thermodynamic model of MU71 describes the time evolution of sea ice temperature profile in the vertical axis along its thickness which also evolves in time due to accumulation or ablation caused by energy balance. 
\subsection{Snow Covered Season Model}
\begin{figure}[t]
\centering
\includegraphics[width=2.0in]{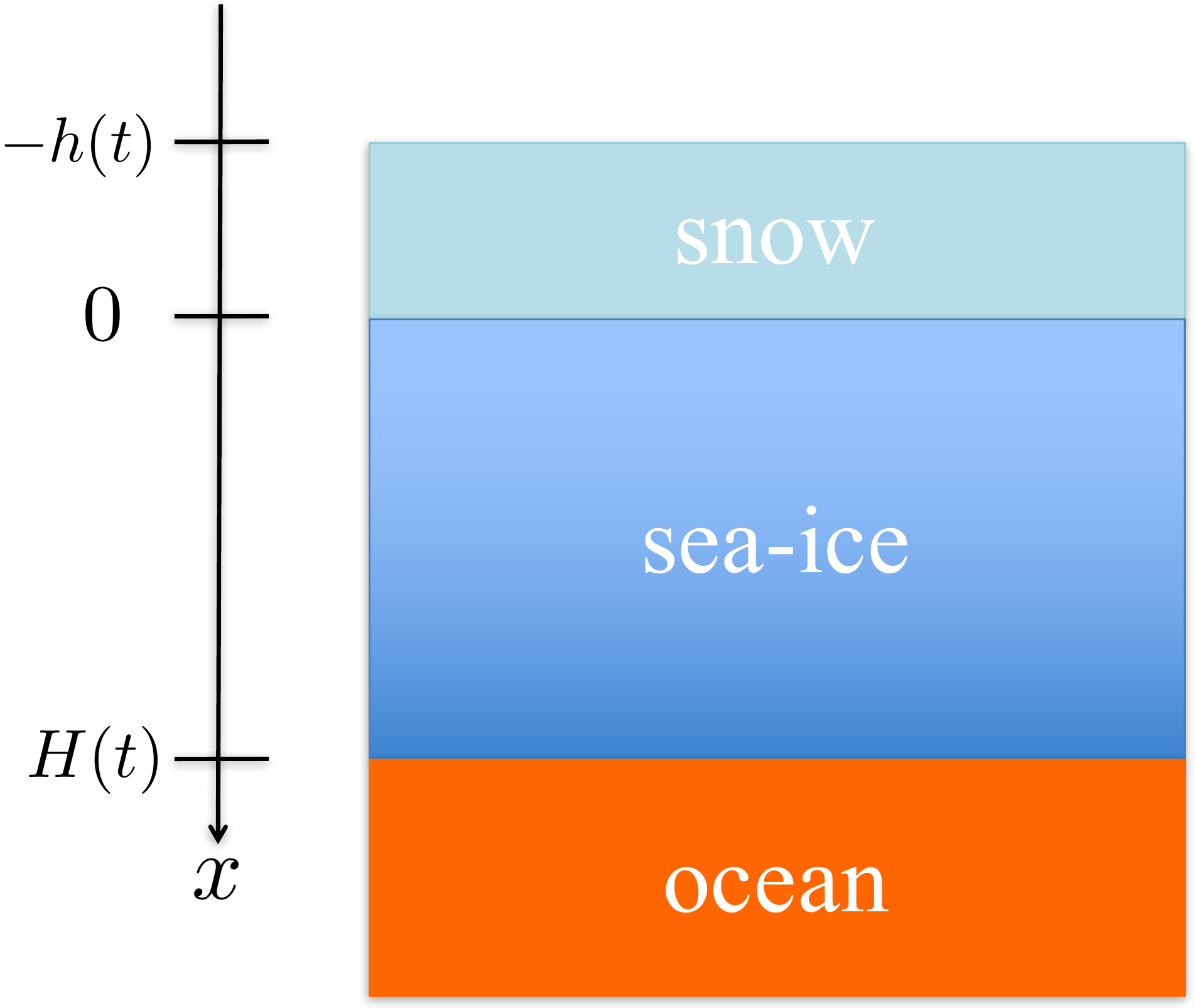}\\
\caption{Schematic of the vertical model of the Arctic sea ice.}
\label{fig:seaice}
\end{figure}

Fig. \ref{fig:seaice} provides a schematic of the Arctic sea ice model. During the seasons other than the summer (July and August), the sea ice is covered by snow, and the surface position of the snow also evolves in time. Let $T_{{\rm s}}(x,t)$, $T_{{\rm i}}(x,t)$ denote the temperature profile of snow and sea ice, and $h(t)$ and $H(t)$ denote the thickness of snow and sea ice. The total incoming heat flux from the atmosphere is denoted by $F_{{\rm a}}$, and the heat flux from the ocean is denoted by $F_{{\rm w}}$. The Arctic sea ice model suggested by MU71 gives the governing equations of a Stefan-type free boundary problem formulated as 
\begin{align}\label{scsys1}
 F_{{\rm a}} - I_0- \sigma (T_{{\rm s}}(-h(t)&,t)+273)^4+ k_{{\rm s}} \frac{\partial T_{{\rm s}}}{\partial x}(-h(t),t)\notag\\
&\hspace{-16mm}=
\left\{
\begin{array}{ll}
0, \quad &{\rm if}\quad  T_{{\rm s}}(-h(t),t) < T_{{\rm m}1}, \\
- q \dot{h}(t) , \quad &{\rm if} \quad T_{{\rm s}}(-h(t),t) = T_{{\rm m}1},
\end{array}
\right.\\
\label{scsys2} \rho_{{\rm s}} c_0 \frac{\partial T_{{\rm s}}}{\partial t}  (x,t)=& k_{{\rm s}}  \frac{\partial^2 T_{{\rm s}}}{\partial x^2}(x,t), \hspace{0mm} \forall x \in (-h(t),0),\\
\label{scsys3} T_{{\rm s}}(0,t) =& T_{{\rm i}}(0,t), \\
\label{scsys3-2} k_{{\rm s}} \frac{\partial T_{{\rm s}}}{\partial x}(0,t) =& k_0 \frac{\partial T_{{\rm i}}}{\partial x}(0,t),\\
\label{scsys4}\hspace{-1mm} \rho c_{{\rm i}} (T_{{\rm i}}, S)\frac{\partial T_{{\rm i}} }{\partial t}(x,t)=& k_{{\rm i}}(T_{{\rm i}}, S) \frac{\partial^2 T_{{\rm i}} }{\partial x^2}(x,t)\notag\\
&+ I_0 \kappa_{{\rm i}} e^{-\kappa_{{\rm i}} x}, \hspace{2mm} \forall x \in (0,H(t)),\\
\label{scsys5}T_{{\rm i}}(H(t),t) =& T_{{\rm m}2}, \\
\label{scsys6}q\dot{H}(t) =& k_{{\rm i}} \frac{\partial T_{{\rm i}}}{\partial x}(H(t),t)-F_{{\rm w}},
\end{align}
where $I_0$, $\sigma$, $k_{{\rm s}}$, $\rho_{{\rm s}}$, $c_0$, $k_0$, $\rho$, $T_{{\rm m}1}$, and $T_{{\rm m}2}$ are solar radiation penetrating the ice, Stefan-Bolzman constant, thermal conductivity of snow, density of snow, heat capacity of pure ice, thermal conductivity of pure ice, density of pure ice, melting point of surface snow, and melting point of bottom sea ice.  The total heat flux from the air $F_{{\rm a}}$ includes the following terms
\begin{align}\label{Fair}
F_{{\rm a}} = (1-\alpha ) F_{{\rm r}}  + F_{{\rm L}}  + F_{{\rm s}} +F_{{\rm l}}, 
\end{align}
where $F_{{\rm r}}$, $F_{{\rm L}}$, $F_{{\rm s}}$, $F_{{\rm l}}$, and $\alpha $ denote the incoming solar short-wave radiation, the long-wave radiation from the atmosphere and clouds, the flux of sensible heat, the latent heat in the adjacent air, and the surface albedo, respectively.  The heat capacity and thermal conductivity of the sea ice are affected by the salinity as 
\begin{align}\label{saltcapa}
c_{{\rm i}}(T_{{\rm i}}, S(x)) =& c_{0} + \gamma_1  \frac{S(x)}{T_{{\rm i}}(x,t)^2}, \\
\label{saltcond}k_{{\rm i}} (T_{{\rm i}}, S(x)) =& k_{0} + \gamma_2 \frac{ S(x)}{T_{{\rm i}}(x,t) },
\end{align}
where $S(x)$ denotes the salinity in the sea ice. $\gamma_1$ and $\gamma_2$ represent the weight parameters. 
The thermodynamic model \eqref{scsys1}-\eqref{scsys6} allows us to predict the future thickness of snow and sea ice $(h(t),H(t))$ and temperature profile $(T_{{\rm s}}, T_{{\rm i}})$ given the accurate initial temperature profile and thickness. However, from the practical point of view, it is not feasible to obtain the complete temperature profile due to a limited number of thermal sensors. To deal with the problem, the estimation algorithm is designed so that the state estimation converges to the actual state starting from an  initial estimate. 

\subsection{Simplification of the Model}\label{Simplification}
Before considering the state estimation design, first we impose a simplification on the system \eqref{scsys1}-\eqref{scsys6}. 
The effect of the salinity profile on the physical parameters is assumed to be sufficiently small so that it can be negligible, i.e.
\begin{align}\label{nosalt}
S(x) = 0.
\end{align}
Therefore, the heat equation of the sea ice temperature \eqref{scsys4} is rewritten as
\begin{align} \label{simpsys4} 
\frac{\partial T_{{\rm i}} }{\partial t}(x,t) =& D_{{\rm i}} \frac{\partial^2 T_{{\rm i}} }{\partial x^2}(x,t) + \bar{I}_0 \kappa_{{\rm i}} e^{-\kappa_{{\rm i}} x}, \hspace{1mm} \forall x  \in (0,H(t)),
\end{align}
where the diffusion coefficient is defined as $D_{{\rm i}} = k_0/\rho c_0$. Next, we impose the following assumptions. 
\begin{assum} \label{ass:Ht} 
The sea ice thickness $H(t)$ is positive and upper bounded, i.e. there exists $\bar{H}>0$ such that 
\begin{align}
0< H(t) < \bar{H}, \quad \forall t>0. 
\end{align}
\end{assum} 
\begin{assum}  \label{ass:dotH} 
The velocity of the thickness $\dot H(t)$ is bounded, i.e., there exists a positive constant $M>0$ such that 
\begin{align} 
 | \dot{H}(t) | < M ,  \quad \forall t>0. 
\end{align} 
\end{assum} 
We formulate the observer structure for sea ice temperature estimation based on the simplified sea ice model composed of \eqref{simpsys4} and \eqref{scsys5}.

\section{State Estimation Design}\label{sec:3}
In this section, we derive the estimation algorithm utilizing some available measurements and show the exponential convergence of the designed estimation to the simplified sea ice temperature. The ice thickness and surface temperature are measured in several studies \cite{Hall2004,Kwok2009,Laxon2013,Rothrock2008}. 
\subsection{Observer Structure}
Suppose that the sea ice thickness and the ice surface temperature are obtained as available measurements ${\mathcal Y}_1(t)$ and ${\mathcal Y}_2(t)$, i.e.
\begin{align}\label{meas1}
 {\mathcal Y}_{1}(t) =& H(t), \\
   \label{meas2}{\mathcal Y}_{2}(t) =& T_{{\rm i}}(0,t) . 
 \end{align}
The state estimate $\hat{T}_{{\rm i}}$ of the sea ice temperature is governed by a copy of the plant \eqref{simpsys4} and \eqref{scsys5}-\eqref{scsys6} plus the error injection of $H(t)$, namely, as follows: 
\begin{align}
 \label{obsys1}\frac{\partial \hat{T}_{{\rm i}}}{\partial t}(x,t) =& D_{{\rm i}} \frac{\partial^2 \hat{T}_{{\rm i}}}{\partial x^2}(x,t) + \bar{I}_0 \kappa_{{\rm i}} e^{-\kappa_{{\rm i}} x}\notag\\
& - p_1(x,t)\left({\mathcal Y}_1(t) - \hat H(t) \right), \hspace{0mm} \forall  x \in (0, H(t)) \\
\label{obsys2}
 \hat{T}_{{\rm i}}(0,t) =& {\mathcal Y}_{2}(t) - p_2(t) \left({\mathcal Y}_1(t) - \hat H(t) \right), \\
 \label{obsys3}\hat{T}_{{\rm i}}(H(t),t) =& T_{{\rm m}2} - p_3(t) \left({\mathcal Y}_1(t) - \hat H(t) \right),  \\
\label{obsys4}  \dot{\hat H}(t) =& p_4(t)  \left({\mathcal Y}_1(t) - \hat H(t) \right) \notag\\
 &+  \beta \frac{\partial \hat T_{{\rm i}}}{\partial x}({\mathcal Y}_{1}(t),t)- \frac{F_{{\rm w}}}{q},  
\end{align}
where $\beta := \frac{k_{{\rm i}}}{q}$. 
Next, we define the estimation error states as
\begin{align}
 \tilde T(x,t) :=& -(T_{{\rm i}}(x,t) - \hat{T}_{{\rm i}}(x,t)), \\
  \tilde H(t) :=& H(t) - \hat H (t),
\end{align}
where the negative sign is added to be consistent with the model developed in \cite{Shumon17journal} for the liquid phase. 
Subtraction of the observer system \eqref{obsys1}-\eqref{obsys4} from the system \eqref{simpsys4} and \eqref{scsys5}-\eqref{scsys6} yields the closed-system of estimation error as
\begin{align}
\label{errorsys3}
\frac{\partial \tilde{T} }{\partial t}(x,t) =& D_{{\rm i}} \frac{\partial^2 \tilde T}{\partial x^2}(x,t) - p_1(x,t) \tilde H(t),  \hspace{1mm} \forall x \in (0, H(t))\\
\label{errorsys4} \tilde{T}(0,t) =& - p_2(t) \tilde H(t),\\
\label{errorsys5} \tilde{T}(H(t),t) =& - p_3(t) \tilde H(t),\\
\label{errorsys6} \dot{\tilde{H}}(t) =& - p_4(t)\tilde{H}(t) -\beta \frac{\partial \tilde T}{\partial x}(H(t),t) . 
\end{align}
Our goal is to design the observer gains $p_1(x,t)$, $p_2(t)$, $p_3(t)$, $p_4(t)$ so that the temperature error $\tilde{T}$ converges to zero, which enables the state estimate $\hat{T}_{{\rm i}}$ to track the actual sea ice temperature $T_{{\rm i}}$. The main theorem of this paper is stated as follows. 
\begin{thm}\label{theorem}
Let Assumptions \ref{ass:Ht} and \ref{ass:dotH} hold. Consider the estimation error system \eqref{errorsys3}-\eqref{errorsys6} with the design of the observer gains
\begin{align}\label{p1gain} 
p_1(x,t) =& \frac{c\lambda}{\beta} x \frac{ I_1 \left( z\right)}{z} + \left( \frac{\ep H(t)}{D_{{\rm i}}} - \frac{3}{\beta} \right) \lambda^2 x   \frac{ I_2 \left( z\right)}{z^2}  \notag\\
&  + \frac{\lambda^3}{D_{{\rm i}}\beta} x^3  \frac{ I_3 \left( z\right)}{z^3} , \\
\label{p2gain} p_2(t) =& 0, \\
\label{p3gain} p_3 (t)=&  -\frac{\lambda}{2 \beta} H(t) - \ep , \\
\label{p4gain} p_4(t) =&  c - \frac{\lambda}{2} \left( 1 - \frac{\lambda H(t)^2}{8 D_{{\rm i}}} \right) + \frac{\beta \lambda }{2D_{{\rm i}}} \ep H(t), 
\end{align}
where $\lambda>0$, $c>0$, and $\ep>0$ are positive free parameters, $z$ is defined by  
\begin{align} \label{zdef} 
z: = \sqrt{\frac{\lambda}{D_{{\rm i}}} (H(t)^2-x^2)}, 
 \end{align} 
and $I_j(\cdot)$ denotes the modified Bessel function of the $j$-th kind. Then, there exist positive constants $c^*>0$ and $\tilde M>0$ such that, for all $c>c^*$, the norm 
\begin{align} 
 \Phi(t) : = \int_{0}^{H(t)} \tilde{T}(x,t)^2 dx + \tilde{H}(t)^2 
\end{align} 
satisfies the following exponential decay 
\begin{align} 
{\Phi}(t) \leq \tilde M {\Phi}(0) e^{- \min\{\lambda, c\} t} , 
\end{align} 
namely, the origin of the estimation error system is exponentially stable in the spatial $L_2$ norm.  
\end{thm}

\begin{remark}\emph{
The observer gains \eqref{p1gain}-\eqref{p4gain} include the thickness $H(t)$, so the gains are not precomputed offline, but are easily calculated \emph{online}, along with the state estimation. }
\end{remark}
\begin{remark}\emph{
The estimation of the snow temperature profile is also achievable using the same observer structure as in Theorem \ref{theorem} with measurements of the snow thickness and snow surface temperature. To avoid the lengthy statement, we do not provide it in this paper. }
\end{remark}

\subsection{Gain Derivation via State Transformation} 
The backstepping is a well-known method to design the observer gains for PDEs \cite{miroslav08,miroslav-first,miroslav09,krstic2009delay}. Hereafter, the partial derivative of a variable in $t$ and $x$ are denoted as the variable with a subscript of $t$ and $x$, respectively. For the estimation error system \eqref{errorsys3}--\eqref{errorsys6}, we apply the following invertible transformations: 
\begin{align}\label{trs1}
\tilde T(x,t) =&  w(x,t) - \int_{x}^{H(t)} q(x,y) w(y,t) dy \notag\\
&- \psi (x,H(t)) \tilde H(t), \\
\label{trs2}w(x,t) = & \tilde{T}(x,t) - \int_{x}^{H(t)} r(x,y) \tilde{T}(y,t) dy \notag\\ 
&- \phi (x, H(t)) \tilde H(t), 
\end{align}
which map the estimation error system \eqref{errorsys3}-\eqref{errorsys6} into the following target system: 
\begin{align}\label{tarsys4}
 w_{t}(x,t) =& D_{{\rm i}} w_{xx}(x,t) - \lambda w(x,t) \notag\\
 &- \dot{H}(t)f(x,H(t))\tilde H(t),  \hspace{2mm} \forall x \in (0,H(t))\\
\label{tarsys5}w(0,t) =& 0,\\
\label{tarsys6}w(H(t),t) =&\ep  \tilde {H}(t),\\
\label{tarsys7}\dot{\tilde{H}}(t) =& - c\tilde{H}(t) -\beta w_{x}(H(t),t) ,
\end{align}
where $f(x,H(t))$ is to be determined. Taking the first and second spatial derivatives of the transformation \eqref{trs1}, we get  
\begin{align}\label{trs1x}
\tilde T_x(x,t) =&  w_x(x,t) + q(x,x) w(x,t)  \notag\\
&- \int_{x}^{H(t)} q_x(x,y) w(y,t) dy \notag\\
&- \psi_x (x,H(t)) \tilde H(t), \\
\label{trs1xx}
\tilde T_{xx}(x,t) =&  w_{xx}(x,t) + q(x,x) w_x(x,t) \notag\\
&+\left(q_x(x,x) + \frac{d}{dx} q(x,x) \right) w(x,t)  \notag\\
&- \int_{x}^{H(t)} q_{xx}(x,y) w(y,t) dy \notag\\
&- \psi_{xx}(x,H(t)) \tilde H(t). 
\end{align} 
Next, taking the time derivative of \eqref{trs1} along the solution of the target system \eqref{tarsys4}--\eqref{tarsys7}, we get 
\begin{align} 
\tilde T_t(x,t) 
=& D_{{\rm i}} w_{xx}(x,t) - \lambda w(x,t) - \dot{H}(t)f(x,H(t))\tilde H(t) \notag\\
&- D_{{\rm i}} \int_{x}^{H(t)} q(x,y) w_{yy}(y,t) dy \notag\\
& + \lambda \int_{x}^{H(t)} q(x,y) w(y,t) dy \notag\\
&+ \dot{H}(t) \int_{x}^{H(t)} q(x,y) f(y,H(t)) dy \tilde H(t) \notag\\
&- \dot{H}(t)  q(x,H(t)) w(H(t),t) \notag\\
& + \psi (x,H(t)) ( c\tilde{H}(t) +\beta  w_{x}(H(t),t)) \notag\\
&- \dot{H}(t) \psi_{H}(x,H(t)) \tilde H(t). 
\end{align} 
Taking integration by parts and substituting the boundary condition \eqref{tarsys6}, we get 
\begin{align} 
\tilde T_t(x,t) 
 =& D_{{\rm i}} w_{xx}(x,t) + D_{{\rm i}} q(x,x)w_{x}(x,t) \notag\\
 &  - (\lambda + D_{{\rm i}} q_{y}(x,x) )  w(x,t)  \notag\\
& + ( \beta \psi (x,H(t))- D_{{\rm i}}  q(x,H(t))) w_{x}(H(t),t) \notag\\
& + ( D_{{\rm i}} \ep q_{y}(x,H(t)) + c \psi (x,H(t)) ) \tilde{H}(t)  \notag\\
&+  \int_{x}^{H(t)} (\lambda q(x,y) - D_{{\rm i}} q_{yy}(x,y) ) w(y,t) dy \notag\\
&-  \dot{H}(t) \tilde H(t) \left( \ep q(x,H(t)) + \psi_{H}(x,H(t)) \right. \notag\\
& \left. + f(x,H(t)) - \int_{x}^{H(t)} q(x,y) f(y,H(t)) dy \right)  . \label{trs1t}
\end{align} 
Thus, by \eqref{trs1xx} and \eqref{trs1t}, we have 
\begin{align} 
&\tilde T_t(x,t) - D_{{\rm i}} \tilde T_{xx}(x,t) + p_1(x,t) \tilde H(t) \notag\\
=&  - \left(\lambda + D_{{\rm i}} \frac{d}{dx}q(x,x) \right) w(x,t)  \notag\\
& + ( \beta \psi (x,H(t))- D_{{\rm i}}  q(x,H(t))) w_{x}(H(t),t) \notag\\
& + \left( D_{{\rm i}} \ep q_{y}(x,H(t)) + D_{{\rm i}} \psi_{xx}(x,H(t)) \right. \notag\\
& \left. + c \psi (x,H(t)) + p_1(x,t)\right) \tilde{H}(t)  \notag\\
&+  \int_{x}^{H(t)} (\lambda q(x,y) + D_{{\rm i}} q_{xx}(x,y) - D_{{\rm i}} q_{yy}(x,y) ) w(y,t) dy \notag\\
&-  \dot{H}(t) \tilde H(t) \left( \ep q(x,H(t)) + \psi_{H}(x,H(t)) \right. \notag\\
& \left. + f(x,H(t)) - \int_{x}^{H(t)} q(x,y) f(y,H(t)) dy \right)  . 
\end{align} 
Substituting $x=0$ and $x = H(t)$ into \eqref{trs1}, we get 
\begin{align} 
\tilde T(0,t) + p_2(t) \tilde H(t) =& - \int_{0}^{H(t)} q(0,y) w(y,t) dy\notag\\
 & \hspace{-6mm}+ ( p_2(t) - \psi (0,H(t)) ) \tilde H(t), \\
\tilde T(H(t),t)  + p_3(t) \tilde H(t) =&(\ep  - \psi (H(t),H(t)) \notag\\
&+ p_3(t)) \tilde H(t). 
\end{align} 
Moreover, substituting $x = H(t)$ into \eqref{trs1x} yields 
\begin{align}
&\dot{\tilde{H}}(t) + p_4(t) \tilde H(t) + \beta 	\tilde T_x(H(t),t) \notag\\
=&  \left( p_4(t) + \beta (\ep q(H(t),H(t)) - \psi_x (H(t),H(t))) - c \right) \notag\\
&\tilde H(t). 
\end{align}
Therefore, for the equations \eqref{errorsys3}--\eqref{errorsys6} to hold, the gain kernel functions must satisfy the following conditions: 
\begin{align} \label{qeq1} 
	q_{xx}(x,y) - q_{yy}(x,y) = & - \frac{\lambda}{D_{{\rm i}}} q(x,y) , \\
\label{qeq2} \frac{d}{dx} q(x,x) =& - \frac{\lambda}{2D_{{\rm i}}} , \quad q(0,y) = 0, \\
\label{psieq1} \psi (x,H(t)) =& \frac{D_{{\rm i}}}{\beta} q(x,H(t)), 
\end{align} 
and the observer gains must satisfy 
\begin{align} 
\label{p1cond}  p_1(x,t) = & - D_{{\rm i}} \ep q_{y}(x,H(t)) - D_{{\rm i}} \psi_{xx}(x,H(t)) \notag\\
&- c \psi (x,H(t)), \\
 \label{p2cond}  p_2(t) = & \psi (0,H(t)), \\
 \label{p3cond}   p_3(t)  =& \psi (H(t),H(t)) - \ep, \\
 \label{p4cond}   p_4(t) = & c - \beta (\ep q(H(t),H(t)) - \psi_x (H(t),H(t))), 
\end{align}
and the function $f(x,H(t))$ must satisfy 
\begin{align} \label{fcond} 
 &f(x,H(t)) - \int_{x}^{H(t)} q(x,y) f(y,H(t)) dy \notag\\
 &= 	- \ep q(x,H(t)) - \psi_{H}(x,H(t)). 
\end{align}
The solutions to \eqref{qeq1}--\eqref{psieq1} are uniquely given by 
\begin{align} \label{qxy} 
q(x,y) =& - \frac{\lambda}{D_{{\rm i}}} x \frac{ I_1 \left( \sqrt{\frac{\lambda}{D_{{\rm i}}} (y^2-x^2)} \right)}{\sqrt{\frac{\lambda}{D_{{\rm i}}} (y^2-x^2)}}, \\
\psi(x,H(t)) =&  - \frac{\lambda}{\beta} x \frac{ I_1 \left( \sqrt{\frac{\lambda}{D_{{\rm i}}} (H(t)^2-x^2)} \right)}{\sqrt{\frac{\lambda}{D_{{\rm i}}} (H(t)^2-x^2)}}. 
\end{align}
Then, by taking the derivatives of the obtained gain kernels, we get 
\begin{align} \label{psixx} 
\psi_{xx}(x,H(t)) =&  \frac{\lambda^2}{D_{{\rm i}}\beta} x \left( 3 \frac{ I_2 \left( z\right)}{z^2} -\frac{\lambda}{D_{{\rm i}}} x^2 \frac{ I_3 \left( z\right)}{z^3} \right), \\
q_y(x,H(t)) = &   - \frac{\lambda^2}{D_{{\rm i}}^2}H(t) x \frac{ I_2 \left( z \right)}{ z^2} , \label{qyH} 
\end{align}
where $z$ is defined by \eqref{zdef}. Using \eqref{qxy}--\eqref{qyH}, the conditions \eqref{p1cond}--\eqref{p4cond} are led to the explicit formulations of the observer gains given as  \eqref{p1gain}--\eqref{p4gain}. In the similar manner, the conditions for the gain kernel functions of the inverse transformation \eqref{trs2} are given by 
\begin{align}
\label{req1}  r_{xx}(x,y)- r_{yy}(x,y) =& \fr{\lambda}{D_{{\rm i}}} r(x,y), \\
\label{req2}  \frac{d}{dx} r(x,x) =& \fr{\lambda}{ 2 D_{{\rm i}}} , \quad r(0,y)=0, \\
 \label{phieq1} \phi (x,H(t))  =& \frac{D_{{\rm i}}}{\beta}  r(x,H(t)), 
\end{align}
and, the function $f(x,H(t))$ is obtained by 
\begin{align} \label{fsol} 
	f(x,H(t)) =  r(x,H(t)) p_3(H(t)) + \phi_{H}(x,H(t)). 
\end{align}
The solutions to \eqref{req1}--\eqref{phieq1} are given by 
\begin{align}
\label{rsol}	r(x,y) =& \frac{\lambda}{D_{{\rm i}}} x \frac{ J_1 \left( \sqrt{\frac{\lambda}{D_{{\rm i}}} (y^2-x^2)} \right)}{\sqrt{\frac{\lambda}{D_{{\rm i}}} (y^2-x^2)}}, \\
\label{phisol} \phi(x,H(t)) =& \frac{\lambda}{\beta} x \frac{ J_1 \left( \sqrt{\frac{\lambda}{D_{{\rm i}}} (H(t)^2-x^2)} \right)}{\sqrt{\frac{\lambda}{D_{{\rm i}}} (H(t)^2-x^2)}}. 
\end{align}
Using the solutions \eqref{rsol} and \eqref{phisol}, the function $f(x,H(t))$ is obtained explicitly by \eqref{fsol}. Then, the solution \eqref{fsol} also satisfy the condition \eqref{fcond}. Hence, the transformation from $(\tilde T, \tilde H)$ to $(w, \tilde H)$ is invertible.

\subsection{Stability Analysis}\label{sec:Stability}
In this section, we prove the exponential stability of the origin of the the estimation error system \eqref{errorsys3}-\eqref{errorsys6} in the spatial $L_2$ norm. First, we show the exponential stability of the origin of the target system \eqref{tarsys4}-\eqref{tarsys7}. We consider the following Lyapunov functional 
\begin{align}\label{V1}
V_{1} = \frac{1}{2} ||w||^2, 
\end{align}
where $|| w||$ denotes the spatial $L_2$ norm of $w$, defined by $|| w || := \sqrt{\int_{0}^{H(t)} w(x,t)^2 dx}$. 
Taking the time derivative of \eqref{V1} together with the solution of \eqref{tarsys4}-\eqref{tarsys7} yields 
\begin{align}
\dot{V}_1 =&- D_{{\rm i}} || w_x ||^2 - \lambda || w ||^2 + \ep \tilde H(t) w_{x}(H(t),t) \notag\\
 &- \dot{H}(t)\tilde H(t)  \int_{0}^{H(t)} w(x,t) f(x,H(t))dx, \notag\\
 &+\frac{\dot{H}(t)}{2}\ep^2 \tilde{H}(t)^2. \label{V1dot} 
\end{align} 
Applying Young's and Cauchy-Schwarz inequalities to the last term in \eqref{V1dot} with the help of Assumption \ref{ass:dotH}, we have 
\begin{align} 
&- \dot{H}(t)\tilde H(t)  \int_{0}^{H(t)} w(x,t) f(x,H(t))dx \notag\\
\leq & M |\tilde H(t)| \cdot  \left( \int_{0}^{H(t)} f(x,H(t))^2 dx \right)^{1/2} \cdot ||w || \notag\\
\label{CSineq}  \leq &  \frac{ M^2  \bar f }{2 \lambda} \tilde H(t)^2  + \frac{\lambda}{2} ||w || ^2 , 
\end{align} 
where $\bar f :=\max_{H(t) \in (0, \bar H)} \int_{0}^{ H(t)}  f(x,H(t))^2 dx $. Applying \eqref{CSineq} to \eqref{V1dot}, we obtain the following inequality: 
\begin{align} 
\dot{V}_1 \leq & - D_{{\rm i}} || w_x ||^2 - \frac{\lambda}{2} || w ||^2  + \ep  \tilde H(t) w_{x}(H(t),t) \notag\\
& + \left(\frac{ M^2 \bar f }{2 \lambda} + \frac{M \ep^2}{2} \right) \tilde H(t)^2.  \label{V1dotineq} 
\end{align}
Next, we consider 
\begin{align} \label{V2}
V_{2} = \frac{1}{2} \tilde{H}(t)^2. 
\end{align} 
Taking the time derivative of \eqref{V2} along with \eqref{tarsys7} leads to 
\begin{align} 
\dot{V}_{2} = &- c \tilde H(t)^2 - \beta \tilde H(t)  w_{x}(H(t),t)  . \label{V2dot} 
\end{align} 
Let $V$ be the Lyapunov functional defined by 
\begin{align} \label{Vdef} 
V = V_1 + \frac{\ep}{\beta} V_2 . 
\end{align} 
Combining \eqref{V1dotineq} and \eqref{V2dot}, the time derivative of \eqref{Vdef} is shown to satisfy the following inequality: 
\begin{align}
\dot{V} \leq &- D_{{\rm i}} || w_x ||^2 - \frac{\lambda}{2} || w ||^2  \notag\\
&- \left( \frac{\ep}{\beta} c  -  \frac{ M^2 \bar f }{2 \lambda} - \frac{M \ep^2}{2}\right) \tilde H(t)^2 .  \label{Vdot1} 
\end{align} 
Hence, by choosing the gain parameter $c$ to satisfy 
\begin{align} 
 c >  \frac{\beta  M^2 \bar f}{\ep \lambda}  + M \ep^2, 
 \end{align} 
 the inequality \eqref{Vdot1} is led to 
 \begin{align} 
 \dot{V} \leq&  - \lambda V_1  -  \frac{\ep }{\beta} c  V_2 \notag\\
 \leq& - \min\{\lambda, c\}  V . \label{Vdotineq} 
\end{align} 
Applying comparison principle to the differential inequality \eqref{Vdotineq}, we get 
\begin{align}\label{Vfin}
&V(t) \leq V(0) e^{- \min\{\lambda, c\} t }.
\end{align}
Hence, the target system \eqref{tarsys4}-\eqref{tarsys7} is exponentially stable at the origin in the spatial $L_2$ norm. Due to the invertibility of the transformations \eqref{trs1} and \eqref{trs2}, there exist positive constants $\underline{M}>0$ and $\bar M>0$ such that the following inequalities hold
\begin{align}\label{bound}
&\underline{M} \left( \int_{0}^{H(t)} \tilde{T}(x,t)^2 dx + \tilde{H}(t)^2 \right) \notag\\
& \leq V(t) \leq \bar M  \left( \int_{0}^{H(t)} \tilde{T}(x,t)^2 dx + \tilde{H}(t)^2 \right) . 
\end{align}
Hence, by \eqref{Vfin} and \eqref{bound}, we arrive at
\begin{align}
& \int_{0}^{H(t)} \tilde{T}(x,t)^2 dx + \tilde{H}(t)^2 \notag\\
 \leq & \frac{\bar M}{\underline{M}}\left( \int_{0}^{H(0)} \tilde{T}(x,0)^2 dx + \tilde{H}(0)^2\right) e^{- \min\{\lambda, c\} t}, 
\end{align}
which completes the proof of Theorem \ref{theorem}. Note that the designed backstepping observer achieves faster convergence with a possibility of causing overshoot since the overshoot coefficient $\bar M/ \underline{M}$ is a monotonically increasing function in the observer gains' parameters $(\lambda, c)$.

\section{Numerical Simulation}\label{sec:4}

While we have focused on the simplified PDE \eqref{simpsys4} to derive a rigorous proof of the proposed state estimation design \eqref{obsys1}-\eqref{obsys4} with observer gains given by \eqref{p1gain}--\eqref{p4gain}, simulation studies are performed by applying  the estimation design to the original thermodynamic model \eqref{scsys1}-\eqref{scsys6} including salinity. 

\begin{figure*}[t]
\begin{center}
\subfloat[Dynamic behavior of snow and sea ice. ]{\includegraphics[width=3.2in]{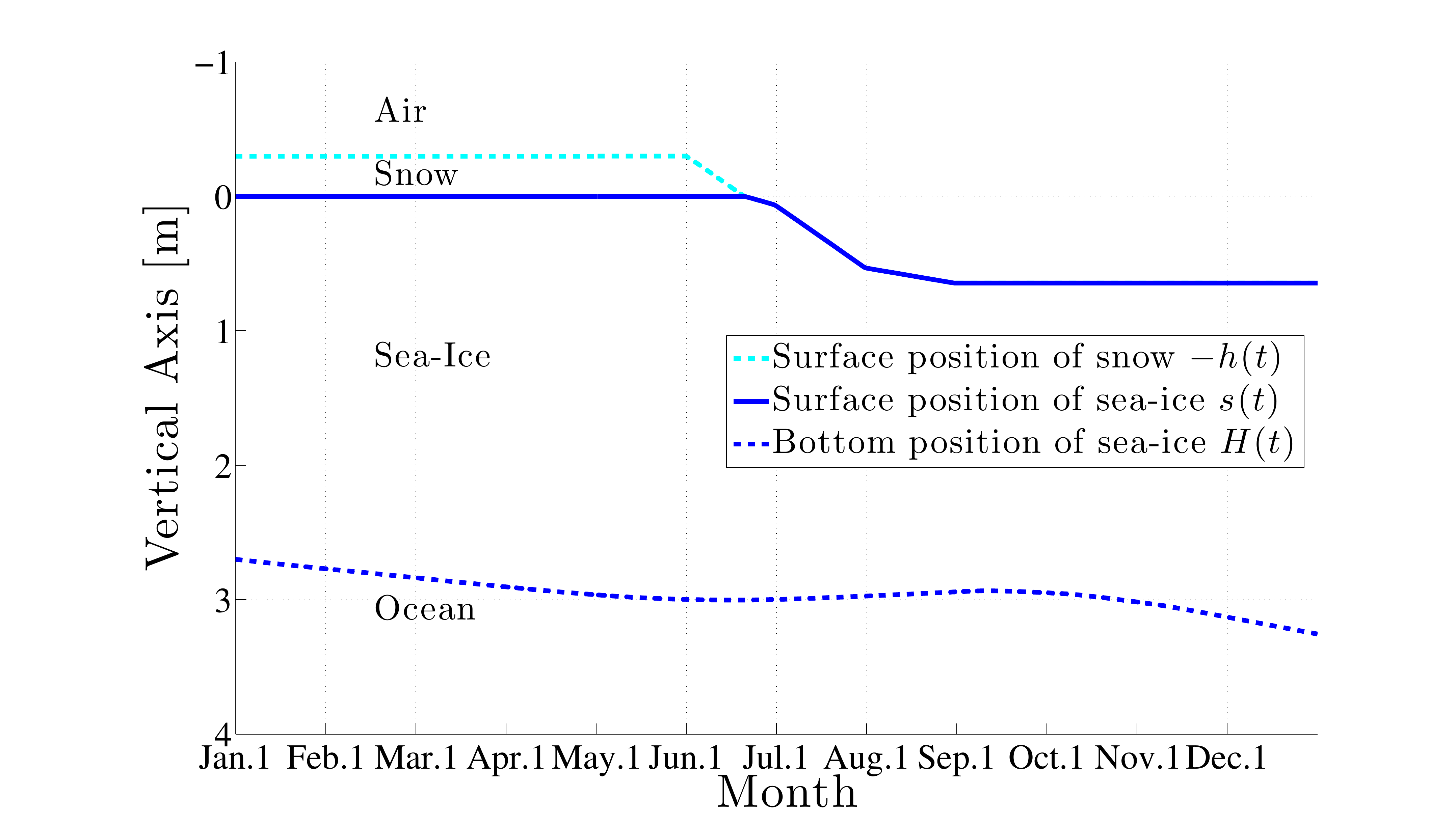}}
\subfloat[Time evolution of temperature profile in sea ice]{\includegraphics[width=3.2in]{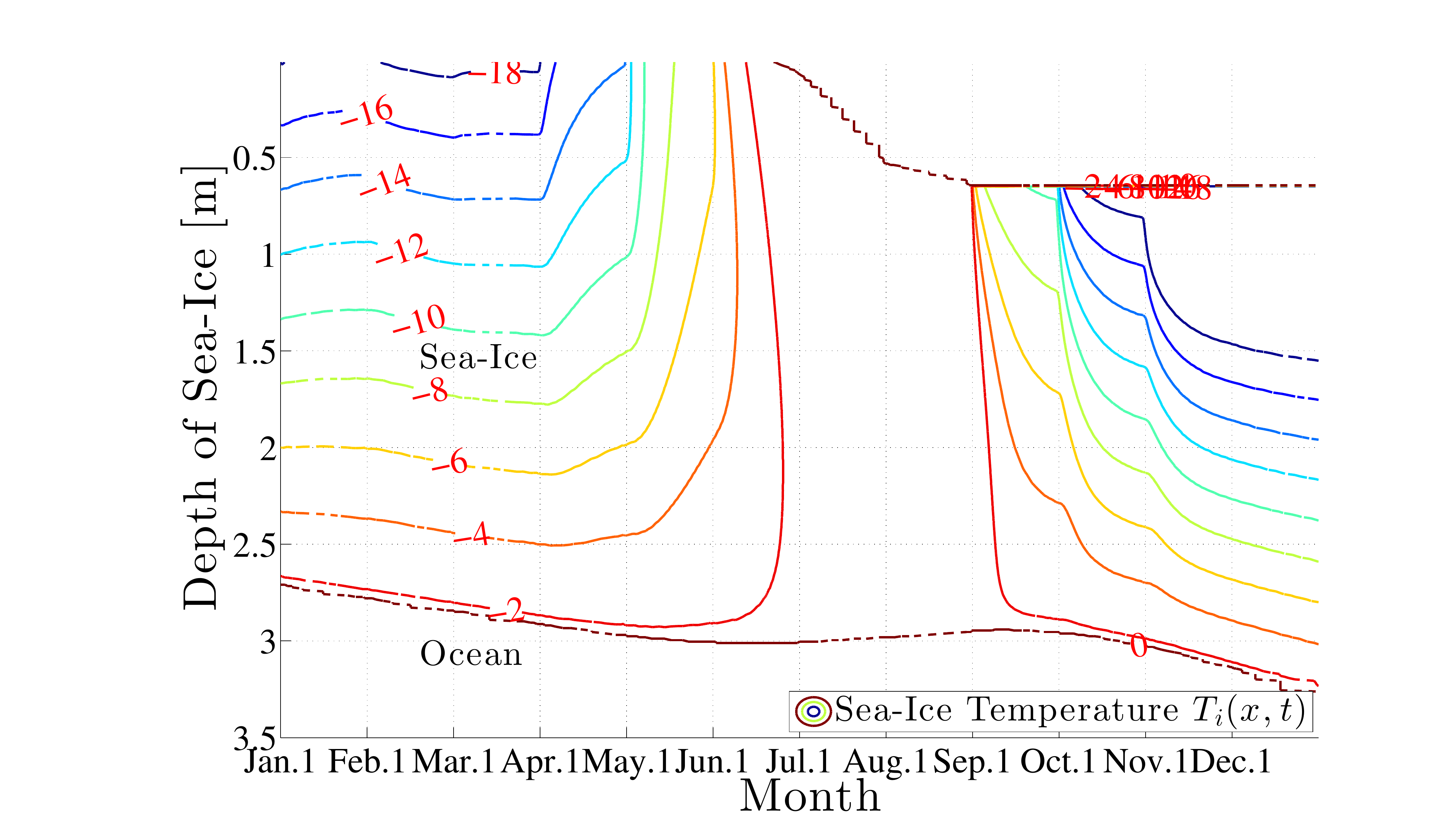}}
\caption{Simulation tests of the plant \eqref{scsys1}--\eqref{scsys6} with input parameters on annual cycle. Both (a) and (b) are in good agreement with the simulation results in \cite{maykut1971}.}
\label{fig:1}
\end{center}
\end{figure*}

\subsection{Input Parameters}

\begin{table}[htb]
\begin{center}
\caption{Average monthly values for the energy fluxes from the atmosphere}

    \begin{tabular}{|c |c c c c c|}
    \hline
     $\textbf{Symbol}$&$F_{{\rm r}}$  & $F_{{\rm L}}$ & $F_{{\rm s}}$ & $ F_{{\rm l}}$ & $\alpha$ \\ \hline
    $\textbf{Unit}$ & W/m$^2$ & W/m$^2$ & W/m$^2$ & W/m$^2$ & \\ \hline
    $\textbf{Jan.}$ & 0 & 168 & 19.0 & 0 & $\cdots$ \\
    $\textbf{Feb.}$ & 0 &  166 & 12.3 & -0.323 & $\cdots$\\
    $\textbf{Mar.}$ & 30.7 & 166 & 11.6 & -0.484 & 0.83 \\
    $\textbf{Apr.}$ &  160& 187 & 4.68 & -1.45& 0.81 \\ 
    $\textbf{May.}$ & 286 &  244 &-7.26 & -7.43 & 0.82 \\
    $\textbf{Jun.}$ & 310 & 291 & -6.30 & -11.3 &  0.78 \\
    $\textbf{Jul.}$ &  220 & 308 &-4.84 & -10.3 & 0.64 \\
    $\textbf{Aug.}$ & 145 & 302 & -6.46 & -10.7 &  0.69 \\
    $\textbf{Sep.}$ & 59.7 & 266 & -2.74 &  -6.30 & 0.84 \\
    $\textbf{Oct.}$ & 6.46 & 224 & 1.61 & -3.07 & 0.85 \\
    $\textbf{Nov.}$ & 0 & 181 & 9.04 &  -0.161 & $\cdots$ \\
    $\textbf{Dec.}$ & 0 & 176& 12.8 & -0.161 & $\cdots$ \\ \hline
    \end{tabular}
     \label{table:1}
    \end{center}
    \end{table}
    \begin{table}[htb]
    \begin{center}
   
\caption{Physical parameters of snow and sea ice}
    \begin{tabular}{| l | l | l |}
    \hline
     $\textbf{Symbol}$& $\textbf{Unit}$ & $\textbf{Value}$ \\ \hline
     $\rho_{{\rm s}}$ & kg/m$^3$ & 330 \\ 
    $k_{s}$ & W/m/$^{\circ}$C & 0.31\\
    $\rho$ & kg/m$^3$ & 917 \\ 
    $ c_{0}$ & J/kg/$^{\circ}$C & 2110  \\  
    $k_{0}$ & W/m/$^{\circ}$C & 2.034\\ 
    $\gamma_1$ & kJ $^{\circ}$C/kg & 18.0 \\ 
    $\gamma_2$ & W/m &  0.117\\ 
     $I_0$ & W/m$^2$ & 1.59 \\ 
    $\kappa_{{\rm i}}$ & /m & 1.5 \\ 
    $T_{{\rm m}1}$ & $^{\circ}$C & -0.1 \\
    $T_{{\rm m}2}$ & $^{\circ}$C & -1.8 \\  \hline
    \end{tabular}
\label{table:2}
\end{center}
\end{table}

The input parameters are taken from  \cite{maykut1971} in SI units and  Table \ref{table:1} shows the monthly averaged values of heat fluxes coming from the atmosphere for each months. Table \ref{table:2} shows the physical parameters of snow and pure ice (without salinity), weight parameters $\gamma_1$ and $\gamma_2$ of salinity effect on heat capacity and thermal conductivity, radiative energy $I_0$, and melting temperatures at surface $T_{{\rm m}1}$ and bottom $T_{{\rm m}2}$. Following the energy conserving model in \cite{Bitz1999}, the salinity profile is described by
\begin{align}\label{salinity}
S(x) = A\left[1 - {\rm cos} \left\{ \pi \left( \frac{x}{H(t)}\right)^{\frac{n}{m+\frac{x}{H(t)}}} \right\}\right], 
\end{align}
where $A = 1.6$, $n = 0.407$, and $m = 0.573$.

\subsection{Simulation Test of MU71}
Using the given data, firstly the simulation of \eqref{scsys1}-\eqref{scsys6} is performed and showed in Fig. \ref{fig:1} to recover the evolution of $h(t)$ and $H(t)$ in the annual season as in \cite{maykut1971}. The dynamic behavior of the snow surface and the bottom of sea ice are shown in Fig. \ref{fig:1} (a), and the time evolution of the temperature profile in sea ice is illustrated in Fig. \ref{fig:1} (b). It can be seen that both of Fig. \ref{fig:1} (a) and (b) have a good agreement with the simulation results shown in \cite{maykut1971}. 

\begin{figure*}[t]
\begin{center}
\subfloat[Open-loop estimation, i.e., $p_1(x,t) = 0$ and {\color{white}aaaa} $p_i(t) = 0$ for $i=2,3,4$. ]{\includegraphics[width=3.2in]{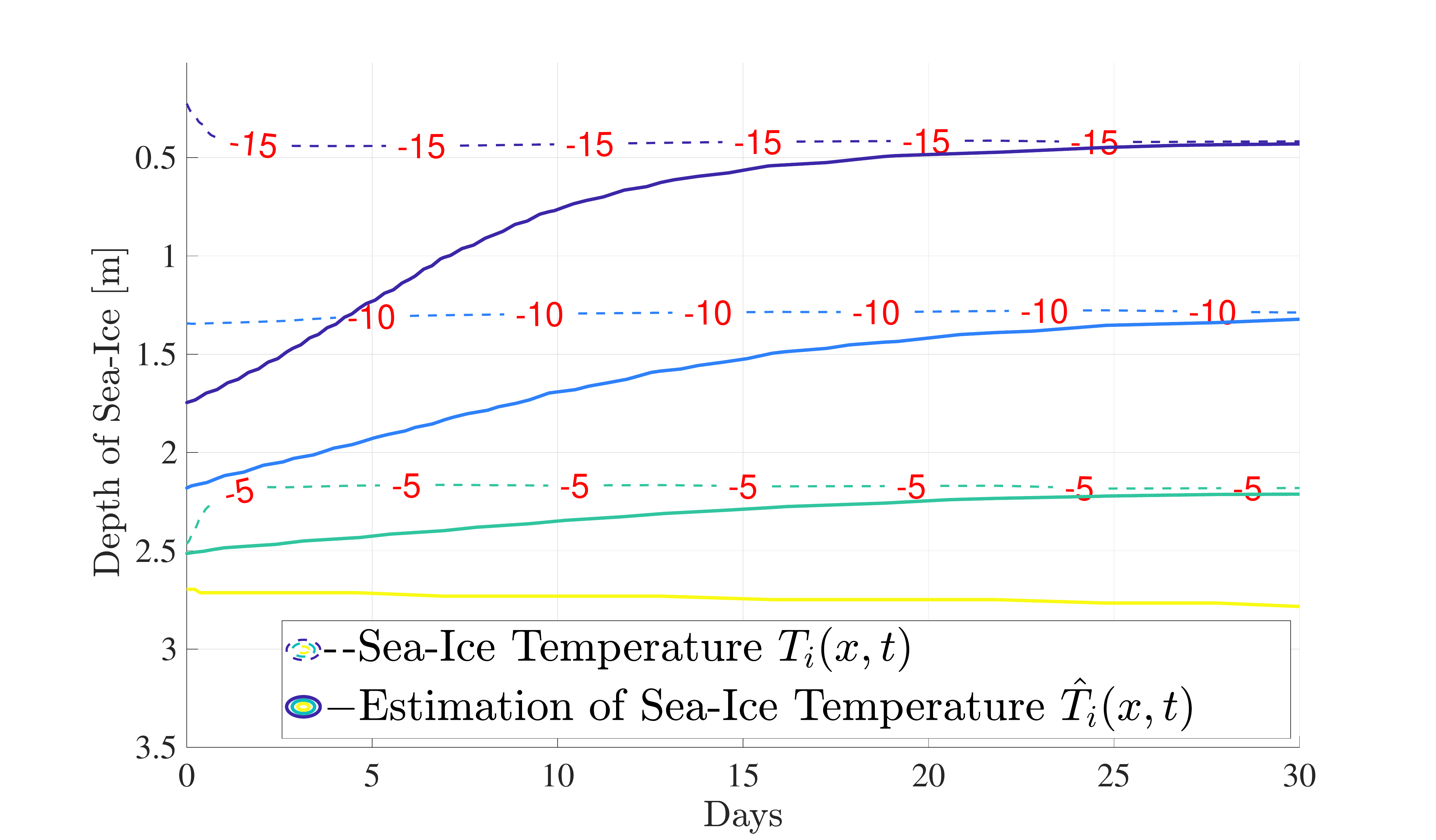}}
\subfloat[The proposed estimation with the observer gains given in \eqref{p1gain}--\eqref{p4gain}. ]{\includegraphics[width=3.2in]{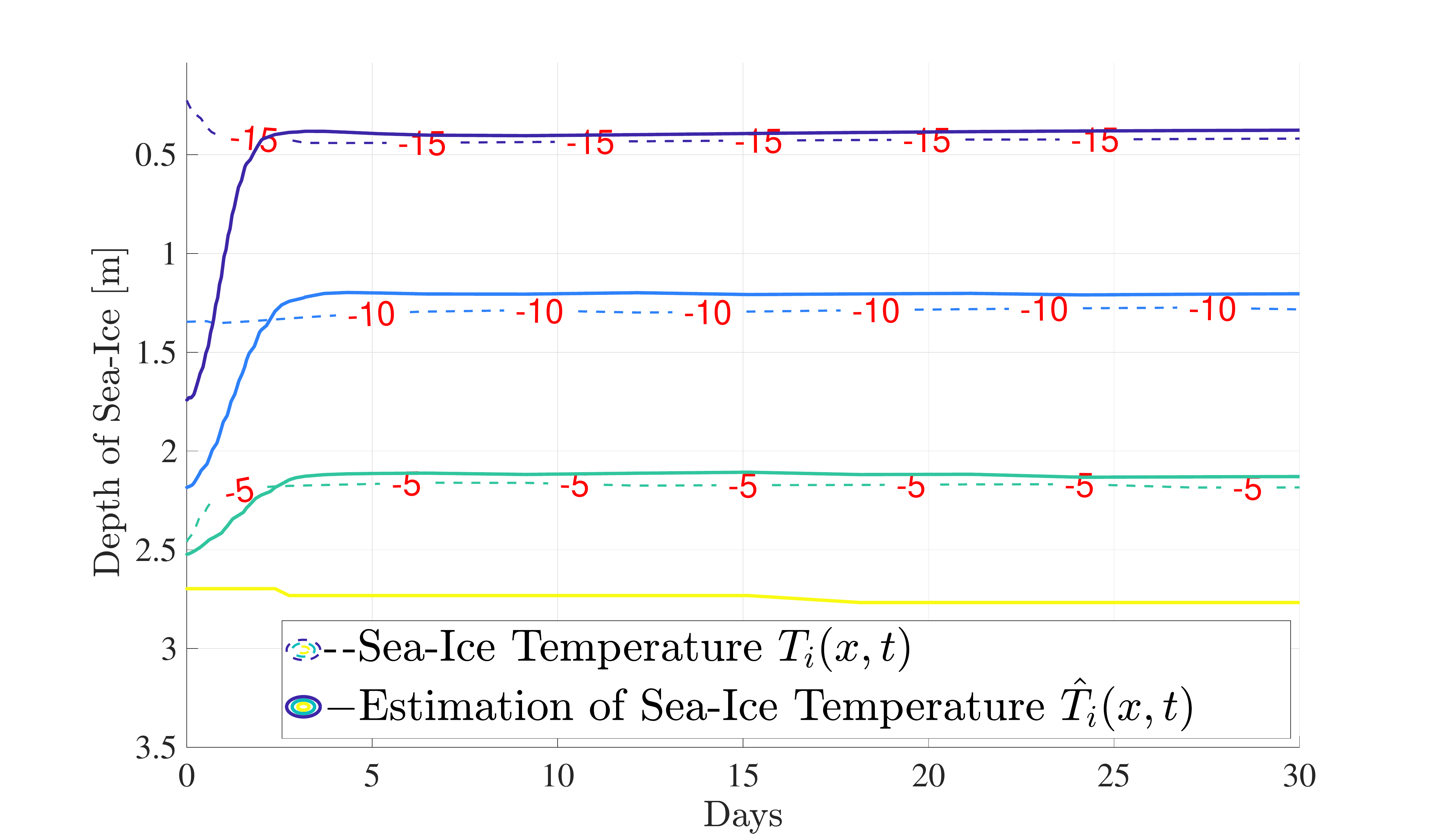}}\\
\subfloat[The proposed estimation with larger value of $\lambda$ {\color{white}aaaa} than the value used in (b). The overshoot is observed. ]{\includegraphics[width=3.2in]{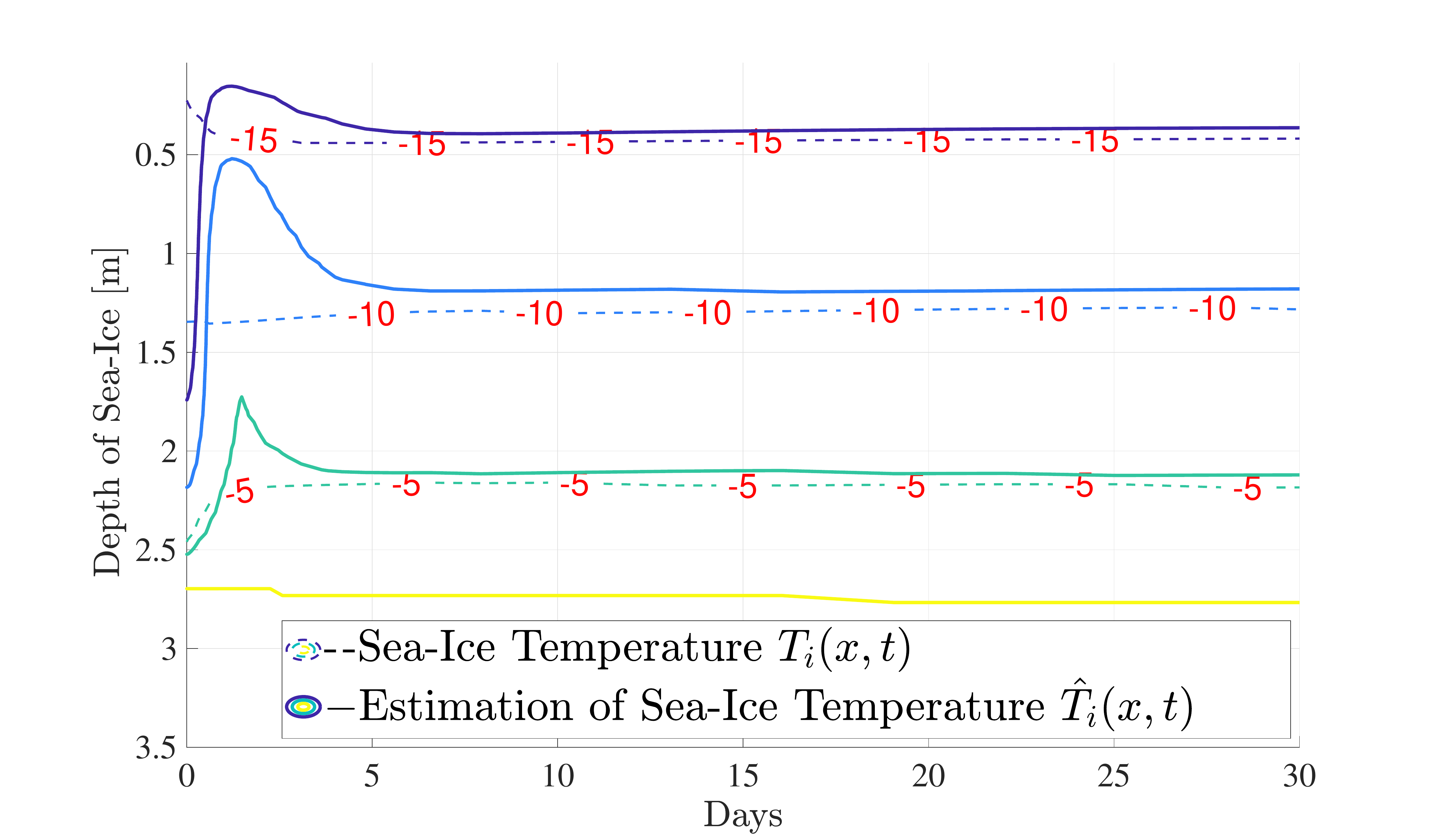}} \subfloat[The proposed estimation with smaller value of $\lambda$ than the value used in (b). The convergence speed gets slower. ]{\includegraphics[width=3.2in]{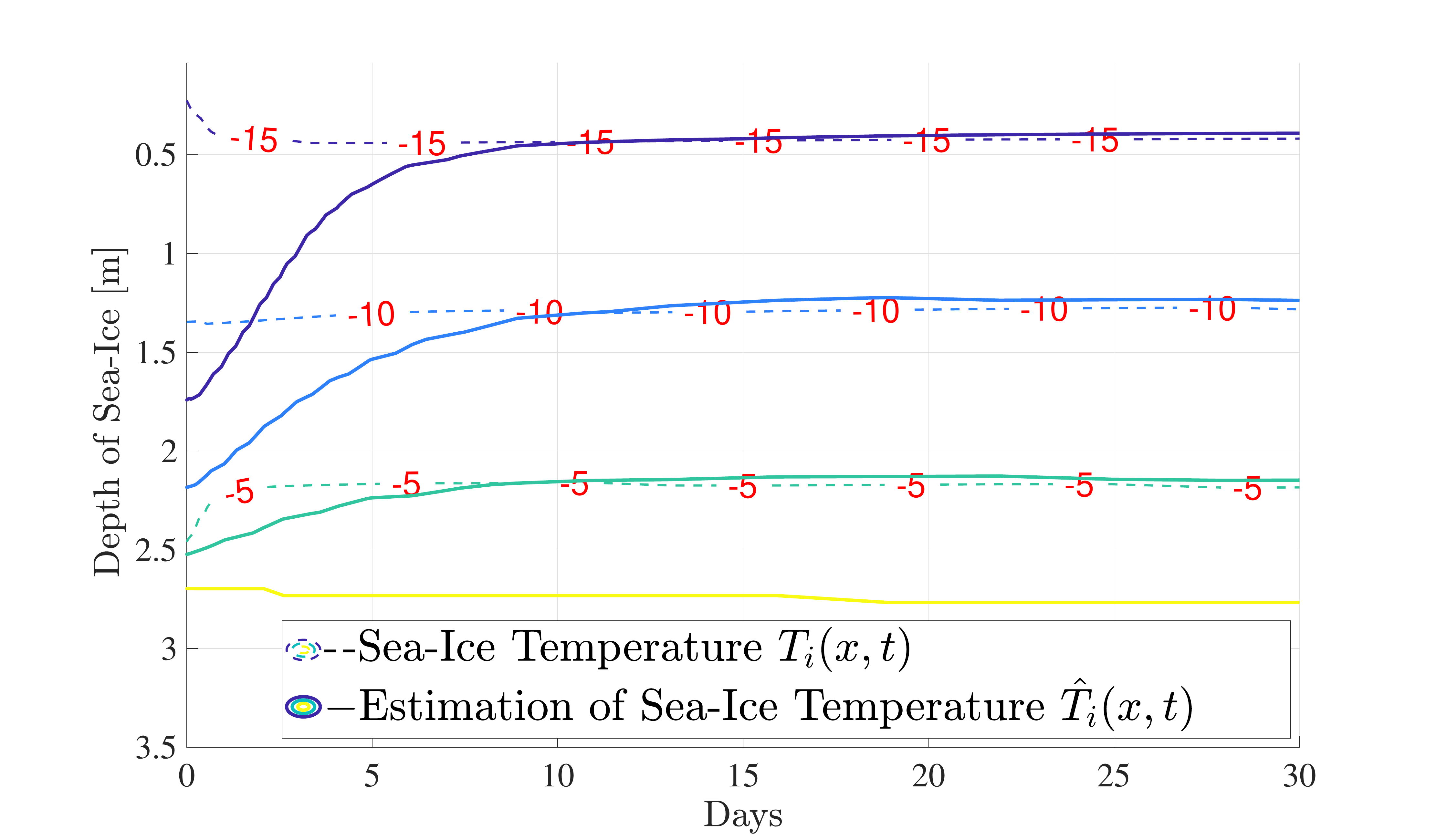}}
\caption{Simulation results of the plant \eqref{scsys1}--\eqref{scsys6} and the estimator \eqref{obsys1}-\eqref{obsys4} using parameters on January. The designed backstepping observer achieves faster convergence to the actual state than the straightforward open-loop estimation.}
\label{fig:2}
\end{center}
\end{figure*}
\begin{figure*}[t]
\begin{center}
\subfloat[Temperature profiles on January 1st.]{\includegraphics[width=3.2in]{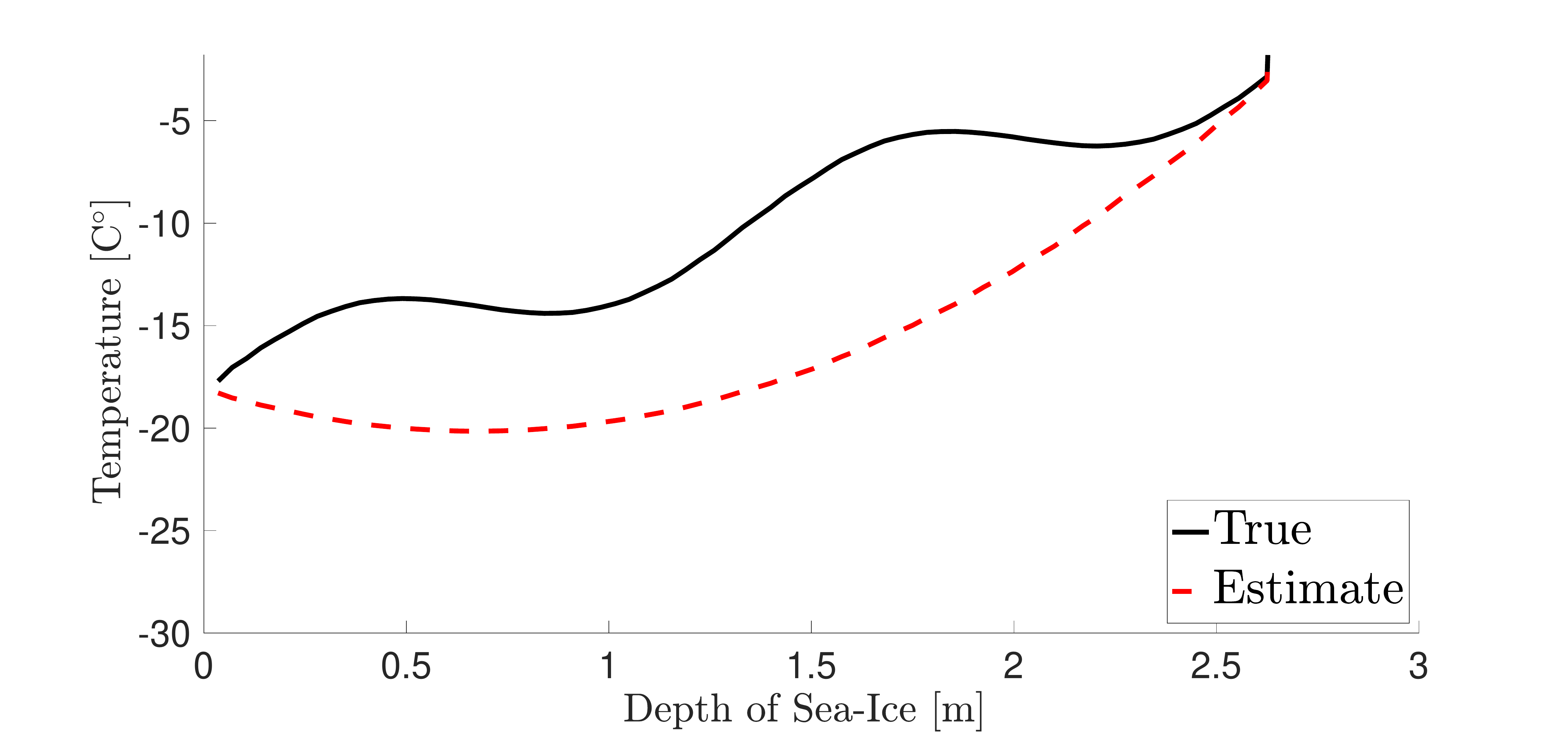}}
\subfloat[Temperature profiles on January 2nd. ]{\includegraphics[width=3.2in]{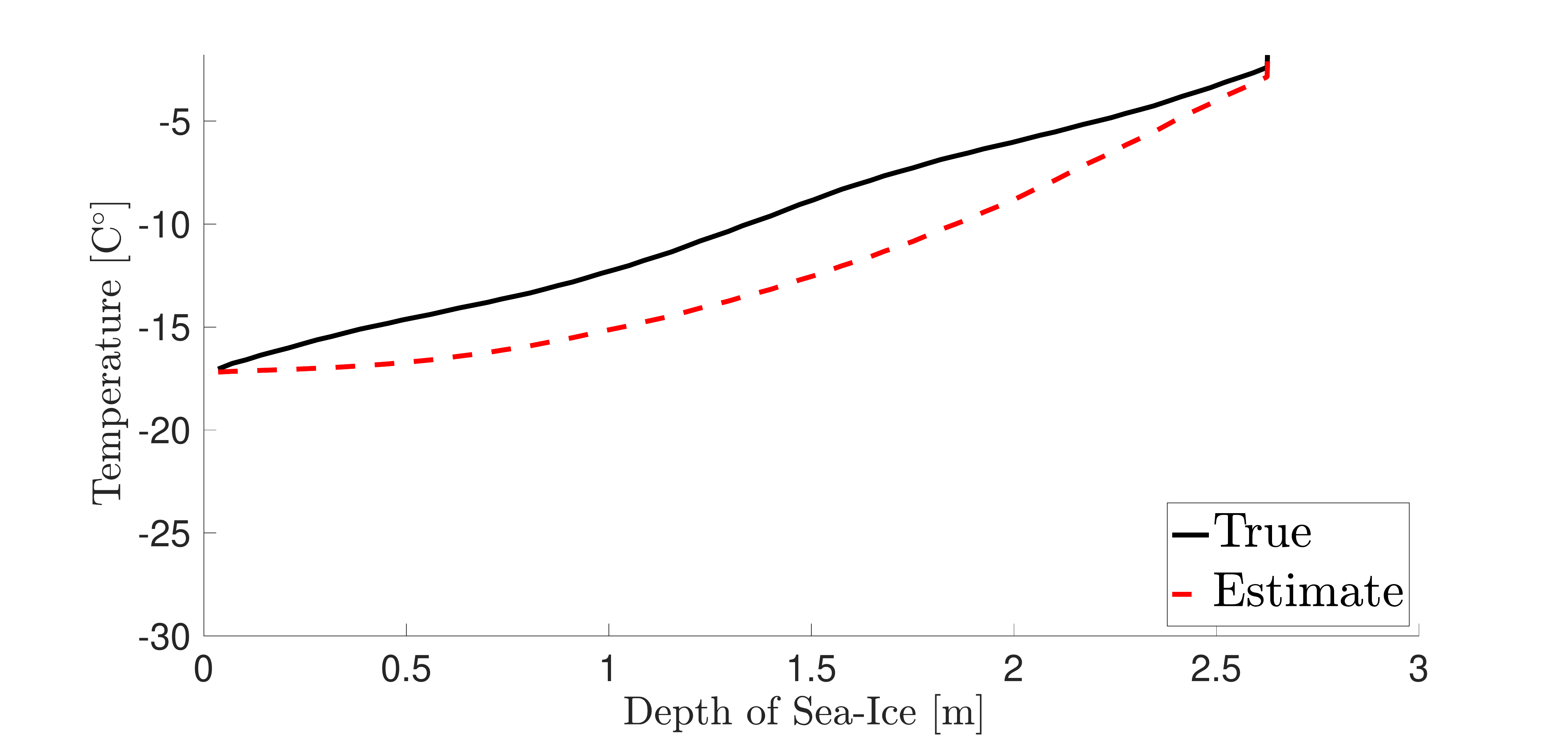}}\\
\subfloat[Temperature profiles on January 3rd.]{\includegraphics[width=3.2in]{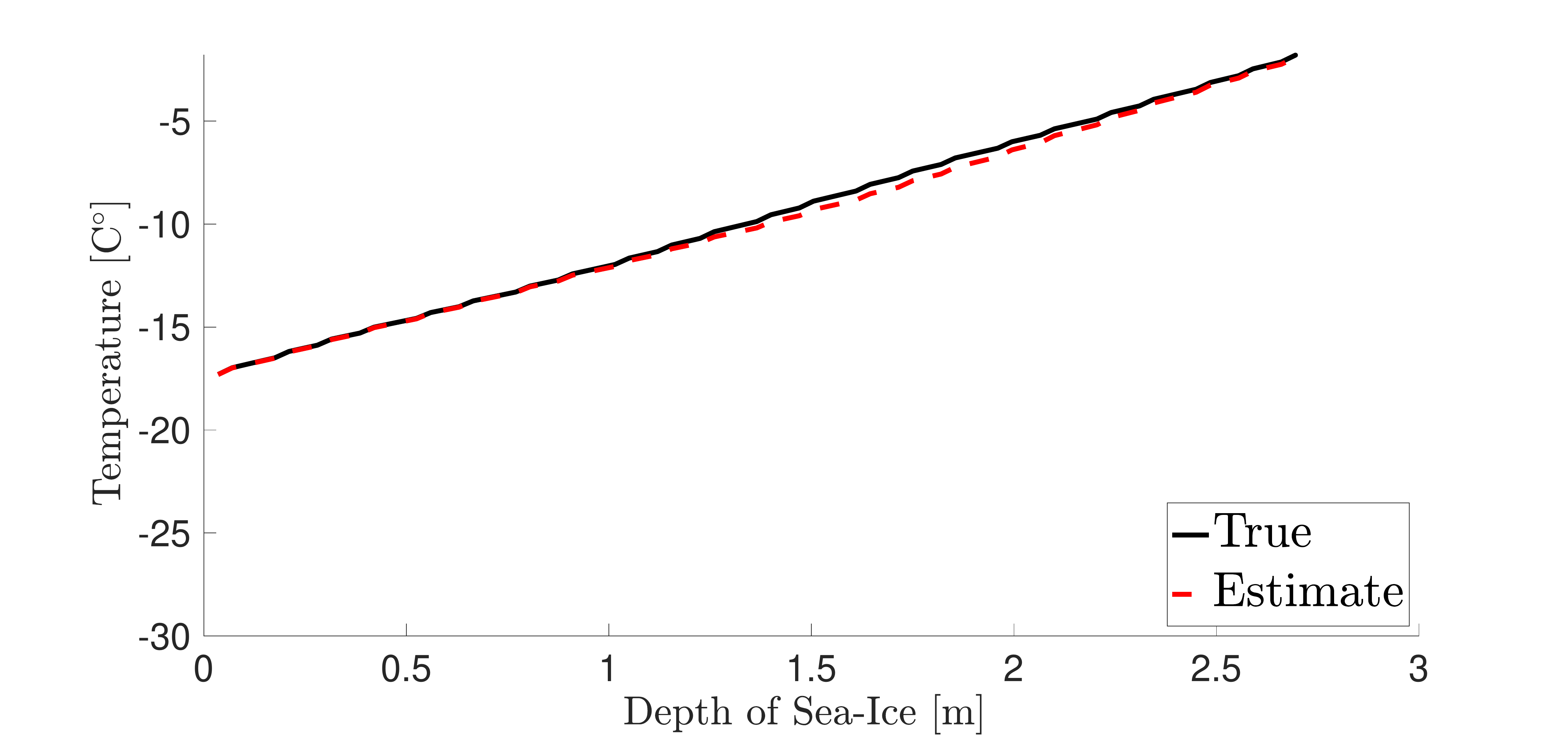}}
\subfloat[The time evolution of $\tilde H(t)$.]{\includegraphics[width=3.2in]{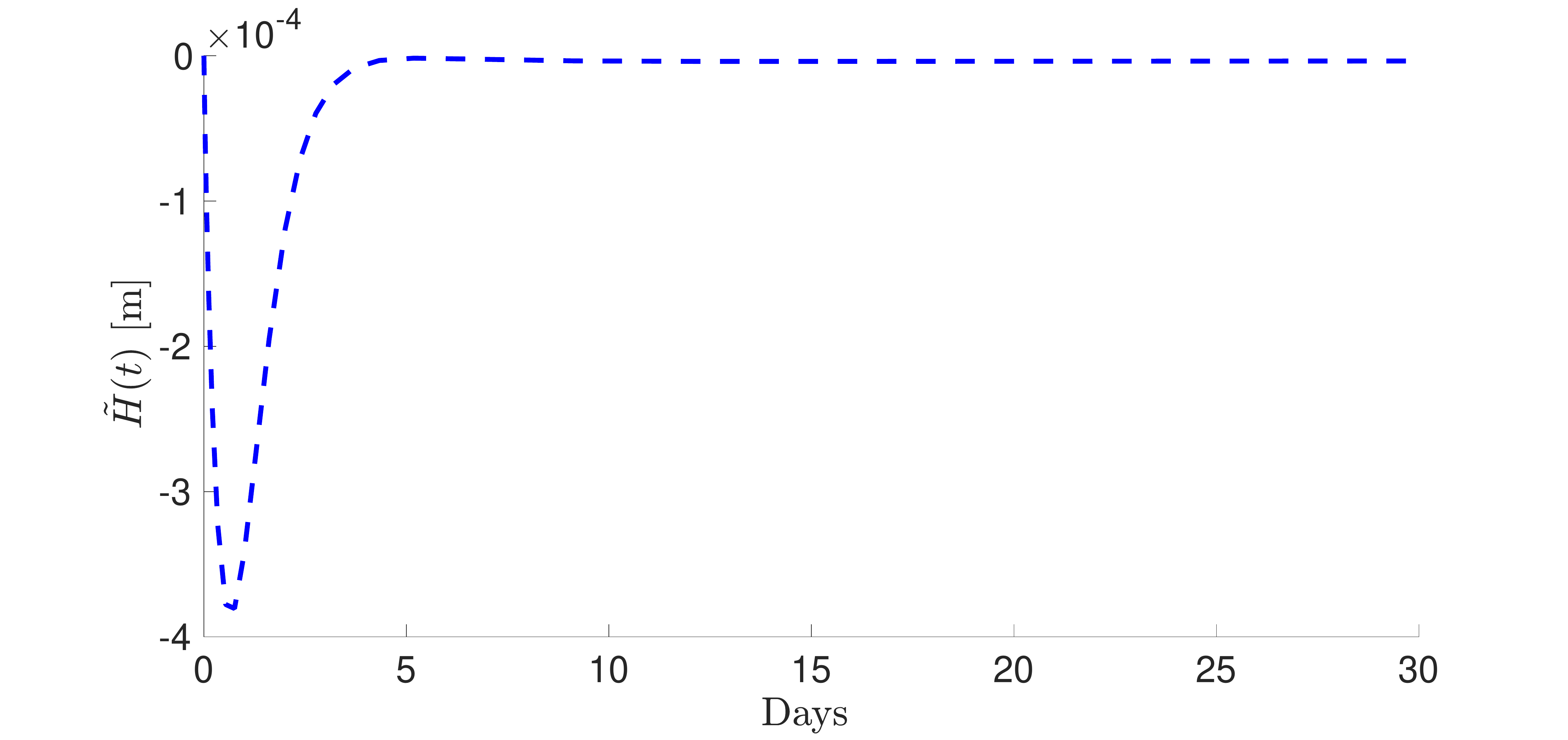}}
\caption{Simulation result of the plant \eqref{scsys1}--\eqref{scsys6} and the estimator \eqref{obsys1}-\eqref{obsys4} with parameters used in Fig. \ref{fig:2} (b). }
\label{fig:3}
\end{center}
\end{figure*}

\subsection{Simulation Results of State Estimation}
The simulation results of temperature estimation $\hat{T}_{{\rm i}}$ computed by \eqref{obsys1}-\eqref{obsys4} along with the available measurements obtained by the online calculation of \eqref{scsys1}-\eqref{scsys6} are shown in Fig. \ref{fig:2}. Here the initial temperature profile is set to be piecewise linear in space plus the sine wave perturbation, formulated as 
\begin{align}
T_{{\rm s}}(x,0) =& \frac{k_{0}(T_{{\rm m}1}-T_{0})}{k_{s} H_{0}} x + T _{0}, \\
T_{{\rm i}}(x,0) =& \frac{T_{{\rm m}1} - T_{0}}{H_{0}} x + T_{0} + a \sin \left( \frac{4 \pi x}{H_0} \right), 
\end{align}
where $T_0 = T_{{\rm i}}(0,0)$ which is obtained by solving fourth order algebraic equation from \eqref{scsys1} and the input data, and $a$ is set as $a = 1$ [C$^\circ$]. The initial temperature estimation is supposed to be a quadratic function in space, given by 
\begin{align}\label{iniestimate}
\hat T_{{\rm i}}(x,0) = \frac{T_{{\rm m}1} - T_{0}}{H_{0}^2 (1-2d)} \left\{(x-dH_0)^2 - (d H_0)^2\right\} + T_{0}, 
\end{align}
 which has a minimum at $x=d H_0$ with a choice of $0\leq d<1/2$. Here we set $d=1/4$. Hence, the initial temperature estimate is lower than the actual temperature. The equation \eqref{iniestimate} satisfies the boundary conditions \eqref{obsys2} and \eqref{obsys3}. The initial state of the estimated ice thickness $\hat H(0)$ is set as that of the true thickness, i.e., $\hat H(0) = H(0)$, which is feasible because the thickness is actually measured. 
 
         Fig. \ref{fig:2} illustrates the contour plots of the simulation results of $T_{{\rm i}}(x,t)$ and $\hat{T}_{{\rm i}}(x,t)$ for open-loop estimation in (a) by setting all the observer gain to be zero,  and the proposed estimation in (b)--(d) with observer gains \eqref{p1gain}--\eqref{p4gain}, respectively, by using input data on January. For the proposed estimation, we fix the parameters of $c = $3.0 $\times$ 10$^{-5}$ and $\ep = $ 1 for all (b)--(d), and use the parameter of $\lambda = $5.0 $\times$ 10$^{-6}$ in (b), $\lambda = $1.0 $\times$ 10$^{-5}$ in (c), and $\lambda = $5.0 $\times$ 10$^{-7}$ in (d). The figures show that the bacsktepping observer gain makes the convergence speed of the estimation to the actual value approximately 5 to 10 times faster at every point in sea ice. As seen in (b)--(d), while the larger choice of $\lambda$ makes the convergence speed faster, it causes more overshoot beyond the actual temperature. Hence, the tradeoff between the convergence speed and overshoot can be handled by tuning the gain parameter $\lambda$ appropriately, thereby the parameters used in (b) achieves the desired performance. The overshoot behavior is noted at the end of Section \ref{sec:Stability} from a theoretical perspective. Consequently, the stability properties stated in Theorem \ref{theorem} for simplified model can be observed in numerical results of the proposed estimation applied to the original model \eqref{scsys1}-\eqref{scsys6}. 
         
         To visualize the convergence of the estimated temperature profile used in (b) more clearly, Fig. \ref{fig:3} illustrate the profiles of both true temperature (black solid) and estimated temperature (red dash) on January 1st to 3rd in (a)--(c), respectively. We can observe that the estimated temperature profile becomes almost same as the true temperature profile on January 3rd, which is two days after the estimation algorithm runs. Moreover, Fig. \ref{fig:3} (d) depicts the time evolution of $\tilde H(t)$ which is an estimation error of the ice's thickness. We can observe that the error is ``enlarged" from $\tilde H(0) = 0$ due to the error of temperature profile, and returns to zero after the temperature profile gets almost same on January 3rd, from which the necessity of the estimator of the ice's thickness is ensured while the thickness is actually measured.

\section{Conclusion and Future Work}\label{sec:5}
In this paper, we develop the estimation algorithm for temperature profile in the Arctic sea ice via backstepping observer design. The observer gains are derived so that the convergence of the state estimate to the actual state is guaranteed theoretically for simplified model. Numerical simulation is employed to test the performance of the observer design with the original thermodynamic model, which illustrates ten times faster convergence of state estimation to the actual temperature than the straightforward open-loop estimation. 

While we have assumed the \emph{online} availability of the measurements of the ice's surface temperature and thickness in this paper, these data acquired by satellites typically accompany a time-delay due to the communication transmitted throughout outer space. Such a time-delay could be compensated by extending the method developed in \cite{koga_2019delay} for control design under actuator delay to the estimator design under sensor delay following the procedure in \cite{miroslav-first, krstic2009delay}. In addition, the physical parameters used in this paper are uncertain variables in practice, and hence it is significant to design a simultaneous state and parameter estimation algorithm. Such an adaptive estimation for parabolic PDEs have been developed in literature such as \cite{Scott14} using reduced-order model via Pade-approximation and \cite{benosmanbook} using data-driven extremum seeking as an iterative learning method. These two extensions will be considered as our future work. 

\section*{Acknowledgement}
The authors would like to thank I. Eisenman for suggesting the sea ice model we considered throughout this paper and helpful discussions regarding the simulation results. The authors would like to thank I. Fenty for enriching our knowledge on recent trend of sea ice state estimation developed in NASA Jet Propulsion Laboratory.


\bibliographystyle{unsrt}

\end{document}